\documentstyle[amsmath,amssymb]{article}

\begin{document}

 \def\GL{{\rm GL}}
 \def\PGL{{\rm PGL}}
 \def\C{{\Bbb C}}
 \def\rr{\rightrightarrows}
 \def\rrr{\Rrightarrow}
  \def\GA{{\rm GA}}
  \def\null{{\rm null}}

  \def\cP{{\cal P}}
   \def\cQ{{\cal Q}}
  \def\cR{{\cal R}}
   \def\cT{{\cal T}}
   \def\cO{{\cal O}}
   \def\cA{{\cal A}}
    \def\cB{{\cal B}}
    \def\cL{{\cal L}}
\def\cF{{\cal F}}

\def\Gr{{\rm Gr}}
\def\Clos{{\rm Close} (\Gr_n)}

  \def\Mat{{\rm Mat}}
  \def\PMat{{\Bbb P}{\rm Mat}}

  \def\phi{\varphi}
\def\epsilon{\varepsilon}
\def\le{\leqslant}
\def\ge{\geqslant}

  \def\SU{{\frak S\frak U}}
  \def\sp{{\frak sp}}

  \def\Hinge{{\rm Hinge}}

 \newcommand{\ord}{\mathop{\rm ord}\nolimits}
\newcommand{\Ker}{\mathop{\rm Ker}\nolimits}
  \renewcommand{\Im}{\mathop{\rm Im}\nolimits}
    \newcommand{\Dom}{\mathop{\rm Dom}\nolimits}
     \newcommand{\Indef}{\mathop{\rm Indef}\nolimits}
     \newcommand{\rk}{\mathop{\rm rk}\nolimits}
      \newcommand{\Mor}{\mathop{\rm Mor}\nolimits}
        \newcommand{\End}{\mathop{\rm End}\nolimits}
          \newcommand{\Aut}{\mathop{\rm Aut}\nolimits}
              \newcommand{\graph}{\mathop{\rm graph}\nolimits}
\def\change{{\rm cha}}
\def\ext{{\rm ext}}
\def\inte{{\rm int}}

 \def\konets{\hfill $\square$}

 \newcounter{sec}
 \renewcommand{\theequation}{\arabic{sec}.\arabic{equation}}

 \begin{center}
{\large\bf Geometry of   $\GL_n(\C)$ on infinity: hinges,
 complete collineations, projective
compactifications, and universal boundary

\medskip

\sc Yurii A. Neretin\footnote%
{supported by the
grants RFBR 98-01-00303 and
NWO 047-008-009}}

 \end{center}

 \bigskip

 Compactifications of semisimple groups
 and of homogeneous spaces arise in  natural ways in
 different branches of mathematics
 (enumerative algebraic geometry, noncommutative
 harmonic analysis, automorphic forms etc.).
  The most important
 construction of this kind is
 the family of objects
  that are called  the {\it Satake boundary}, or
 the {\it De Concini--Procesi boundary}
 or the {\it wonderful compactification},
  see  \cite{Sat},  \cite{Fur},
   \cite{DCP}, \cite{Kno},  see also \cite{Ner5}.
  For the first time, such compactifications
   of the group
  $\PGL(n,\C)$
  (the {\it complete collineations}) and of the
  symmetric space  $\PGL(n,\C)/{\rm PO}(n,\C)$
  (the {\it complete quadrics})
  were discovered by Semple
   \cite{Sem1}, \cite{Sem2}, \cite{Sem3}%
   \footnote{The
   wonderful compactification
   of $\PGL(3,\C)/{\rm PO}(3,\C)$  was constructed by E.Study in 1886.}.

 In the geometry of symmetric spaces and in the analysis
 on symmetric spaces,
  there arise some
 more complicated compactifications as the Karpelevich
 compactification (\cite{Kar}, \cite{Kus3}, \cite{Gui})
 and
 the  Martin compactification (\cite{Dyn}, \cite{Ols},
 \cite{Gui}).

 There exists also a wide theory of algebraic
 equivariant compactifications
 of reductive groups and their homogeneous
 spaces, see \cite{VP}, \cite{LV},
 \cite{Kno}, \cite{Pop}, \cite{BLV}
 (if a reductive group is a torus, then this theory
 becomes the theory of toric varieties).

 This paper has two purposes.
  The first aim is to give an explicit description
  in elementary geometric terms
  of all the algebraic projective compactifications (see below)
  of $\GL_n(\C)$.
  The second aim is to construct a universal object
 (the sea urchin)
  for all the compactifications of this type.

  \smallskip

  {\bf 0.1. Projective compactifications.}
  Consider the group $\GL_n(\C)$ of all complex invertible
  $n\times n$ matrices.
  Consider a polynomial
   (generally speaking, reducible)
   representation $\zeta$ of the group $\GL_n(\C)$ in
   an $N$-dimensional complex linear space $Z$.

   Denote by $\Mat(Z)$ the space of all linear operators
   in $Z$.
   Consider the space $\PMat(Z)$ consisting of nonzero operators
   defined up to a nonzero scalar factor. This space is the complex
   projective space $\C{\Bbb P}^{N^2-1}$. Consider the maps
   $$\GL_n(\C)\stackrel{\zeta}\to \Mat(Z)\to\PMat(Z).$$
   Denote by $[\GL_n]_\zeta$ the  closure of the image
   of $\GL_n(\C)$ in $\PMat(Z)$. The spaces
   $[\GL_n]_\zeta$ are called
    {\it projective compactifications}
   of $\GL_n(\C)$, see \cite{Rash}.

   \smallskip

   {\sc Remark.} Let $\zeta$
   be an irreducible representation
   with a signature $\nu=(\nu_1,\dots,\nu_n)$ (see below 2.12) satisfying
   the condition
   $$\nu_1> \nu_2>\dots> \nu_n.$$
   In this case, the space  $[\GL_n]_\zeta$ is called
   the Semple complete collineation
   space or
   the De Concini--Procesi
   compactification of $\PGL_n(\C)$.   \konets

  \smallskip

  {\bf 0.2. Abstract characterization
   of projective compactifications.}
  Let the group $\GL_n(\C)\times\GL_n(\C)$
  act on a projective algebraic  variety
  $M$. Denote by $\Delta\simeq\GL(n,\C)$
  the diagonal subgroup in  $\GL_n(\C)\times\GL_n(\C)$.
       Let $G$ have
    an open  $\GL_n(\C)\times\GL_n(\C)$-orbit on $M$, and the
   stabilizer of the orbit contain $\Delta$
  (any proper subgroup of $\GL(n,\C)\times\GL_n(\C)$
  containing $\Delta$ is a product of $\Delta$ and a
   subgroup in the center
  $\C^*\times\C^*$ of $\GL(n,\C)\times\GL_n(\C)$).

 The Kambayashi theorem (\cite{Kam}, see also
 the exposition in \cite{CG}, 5.1) implies that
all normal $\GL_n(\C)\times\GL_n(\C)$-varieties
  $M$ satisfying these conditions
are of the form  $[\GL_n]_\zeta$.

  We  exploit only the constructive definition 0.1.

  \smallskip

  {\bf 0.3. Sea Urchin.} The Sea Urchin $\SU_n$ is
  the universal object for all the projective compactifications
  of $\GL_n(\C)$
  in the following  sense.

  \smallskip

  a) For each projective compactification
   $[\GL_n]_\zeta$, there exists
  a canonical $\GL_n(\C)\times\GL_n(\C)$-equivariant surjective  map
   $$\pi_\zeta:\SU_n\to[\GL_n]_\zeta.$$

   b) Denote by  $\dot D_\epsilon$ the punctured disk
   $0<|z|<\epsilon$ in $\C$.
   Consider the germ in 0 of
   an algebraic curve
    $$\dot D_\epsilon\to\GL_n(\C).$$
  Any such germ    has a limit as $z\to 0$
   in the sea urchin $\SU_n$. Two germs $\gamma_1$, $\gamma_2$ have
    the same limit
   in the sea urchin iff for each
   projective compactification
   $[\GL_n(\C)]_\zeta$  the  limits
   $\lim_{z\to 0}\gamma_1(z)$ and
     $\lim_{z\to 0}\gamma_2(z)$  in $[\GL_n(\C)]_\zeta$
     coincide.

      \smallskip

  c) $\pi_\zeta(\lim_{z\to 0} \gamma(z))=\lim_{z\to 0} \pi_\zeta(\gamma(z))$
  for all the germs $\gamma$ and all the $[\GL_n]_\zeta$.

      \smallskip

 The existence of the sea urchin is obvious.
   Strangely enough,
 this object admits an explicit
  parametrization in elementary terms.  Points
  of the sea urchin are
 enumerated by collections of integers $(m_1,\dots, m_n)$
 (defined up to a common factor)
 and some special  collections of subspaces (hinges)
 $P_1,\dots,P_\tau\subset\C^n\oplus\C^n$.

 \smallskip

 {\sc Remark.}
  Obviously,
 there exists a universal object for all the projective
  compactifications
 in the category of compact topological spaces%
 \footnote{Let $M_1,M_2,\dots$ be all the projective compactifications.
 Consider the diagonal embedding $\GL_n\to M_1\times M_2\times\dots$
 (the product is equipped with the topology of pointwise convergence =
 the Tihonov topology).
 The closure of the image  of this
 embedding is the required universal object.}.
 This object
  (for the case of Riemannian noncompact
  symmetric spaces) was investigated by
 Kushner \cite{Kus1}, \cite{Kus2}.  The sea urchin
 is not a compact space in the usual sense,
 hence the sea urchin differs from the Kushner compactification.
 The sea urchin also is simpler.         \konets

   \smallskip

 {\sc Remark.} Obviously, the sea urchin is not a projective variety.
                                              \konets

 \smallskip

 {\bf 0.4. Basic observations.} Denote by $\Gr_n$
 the Grassmannian  of $n$-dimensional
 subspaces in $\C^n\oplus\C^n$. Consider the canonical
 embedding  $\GL_n\to \Gr_n$ taking any operator to its graph.

 Consider the germ in $z=0$
  of an algebraic (or meromorphic)
 map $\gamma:\dot D_\epsilon\to\GL_n(\C)$.
 For any integer $k$, consider the following
 limit in the Grassmannian
 $$R_k=\lim\limits_{z\to 0}
       z^{-k} \gamma(z).$$
 We obtain the bilateral sequence
 \begin{equation}
 \dots, R_{-2},R_{-1},R_0,R_1,R_2,\dots
 \end{equation}
 Consider all $k$ such that $R_{k}$
 is not a sum of a horizontal subspace
 and a vertical subspace, i.e.,
 $$
 R_{k}\ne \Bigl[R_{k}\cap \bigl(\C^n\oplus 0\bigr)\Bigr]
       \Bigl[R_{k}\cap \bigl(0\oplus \C^n\bigr)\Bigr]
.$$
 Thus we obtain some {\it finite} collection of integers
 $$k_1>k_2>\dots> k_\tau.$$
  Then we select the corresponding $R_{k_j}$ from
  the sequence (0.1).
   Thus we obtain the finite
 ($\tau\le n$) family of subspaces in $\C^n\oplus\C^n$
 $$\cR=(R_{k_1},\dots, R_{k_\tau}).$$
 All possible families $\cR$ can  easily be described,
 see the definition of hinges in 2.1.

 Thus, for each meromorphic germ $\gamma(z)$,
 we associate the following data
 \begin{equation}
 (k_1,\dots,k_\tau;\,\, R_{k_1},\dots, R_{k_\tau})
 .\end{equation}

 Our paper contains two  observations.

 \smallskip

 1. For any polynomial representation
    $\zeta$ of $\GL_n(\C)$, the limit of $\zeta(\gamma(z))$
    in  $[\GL_n]_\zeta$ is completely determined
    by the data (0.2)

 \smallskip

 2. The operator
    $$\lim\limits_{z\to 0} \zeta(\gamma(z))$$
    admits a simple explicit description in terms
   of the data (0.2).

  \smallskip

  This gives the explicit description of all the spaces
 $[\GL_n]_\zeta$ and also the description of the sea urchin.

 \smallskip

 {\bf 0.5. Structure of this paper.}
  Section 1 contains preliminaries on
 the category $\GA$ of linear relations
  and on the fundamental representation of $\GA$. These objects
appeared in \cite{Ner1}, the detailed exposition is contained
 in \cite{Ner2} and \cite{Nerb}, 2.5.

  Section 2 contains the preliminaries on the hinges
  and the hinge semigroup.
  It is mainly based on  \cite{Ner3}, except
  Subsections 2.10, 2.11;
  the detailed exposition
   of the work \cite{Ner3} is present in \cite{Ner5}.
   We also explain here relations between the
   hinge semigroup and some well-known constructions.

   Section 3  contains the  construction of
  the sea urchin.

  \smallskip

  {\bf 0.6. Other classical groups and symmetric spaces.}
  In this paper, we consider only the groups
  $\GL_n(\C)$.  The hinge language
is common
   for all the classical groups, all the classical symmetric
   varieties, and for their real forms (i.e., Riemannian
   and pseudo-Riemannian symmetric spaces),
   see \cite{Ner3}, \cite{Ner5}.
   The results of this paper extend to this  general
   situation more or less automatically.

   \smallskip

  {\bf Acknowledgments.}
  I am  grateful to C.De Concini, E.B.Vinberg and V.L.Popov
  for meaningful discussion of this subject.
  I thank the administrators of  the Erwin Schr\"odinger
  Institute
for Mathematical Physics,
where this work was done, for their hospitality.

\smallskip

{\bf Notation.}
We denote by $\C^*$ the multiplicative group of
nonzero complex numbers.

 Let $V$ be an $n$-dimensional
complex linear space. We denote by $\GL_n(\C)$
or $\GL(V)$ the group of all invertible linear operators
in $V\simeq\C^n$.

 By $\PGL_n(\C)$ we denote the quotient group
$\GL_n(\C)/\C^*$ of $\GL_n(\C)$ by the subgroup
$\C^*$ of all scalar matrices.


 For a linear space $Z$, we denote by ${\Bbb P}Z$
 the projective space $(Z\setminus 0)/\C^*$.

 By $\Mat(Z)$ we denote the space of all linear operators
 in a linear space $Z$. By $\PMat(Z)$ we denote the quotient space
 $(\Mat(Z)\setminus 0)/\C^*$.

     \bigskip

 {\large\bf 1. Category $\GA$
 and its fundamental representation}

 \medskip
 \addtocounter{sec}{1}
\setcounter{equation}{0}

 {\bf 1.1. Linear relations.}
 Let $V$, $W$ be finite-dimensional linear spaces over $\C$.
 A {\it linear relation} $P:V\rr W$ is
 a linear subspace in $V\oplus W$.

 \smallskip

 {\sc Example.} Let $A:V\to W$ be a linear operator.
 Its {\it graph}
$\graph (A)$ consists of all vectors
of the form $v\oplus Av\in V\oplus W$.
 Thus $\graph(A)$ is a linear relation $V\rr W$.
 {\it Below we do not distinguish linear operators and their graphs.} \konets

 \smallskip

 Let $V,W,Y$ be linear spaces, and let $P:V\rr W$, $Q:W\rr Y$
 be linear relations. Their product $S=QP$  is the linear
 relation $S:V\rr Y$ consisting of all $v\oplus y\in V\oplus Y$
 such that there exists $w\in W$ satisfying the conditions
 $$ v\oplus w\in P,\qquad w\oplus y\in Q.$$

 \smallskip

 For any  linear relation $P:V\rr W$, we define

  \smallskip

  a) the {\it kernel} $\Ker P\subset V$ is
  $P\cap (V\oplus 0)$;

  \smallskip

  b) the {\it image} $\Im P\subset W$ is the image of the projection
     of $P\subset V\oplus W$ on $W$ along $V$,

  c) the  {\it domain} $\Dom P\subset V$
  is the image of the projection
     of $P\subset V\oplus W$ on $V$ along $W$,

    \smallskip

  d) the {\it indefineteness} $\Indef P\subset W$
  is $P\cap (0\oplus W)$,

   \smallskip

   e) the {\it dimension} $\dim P$ is the dimension of $P$,

   \smallskip

   f) the {\it rank}
 \begin{multline*}
 \rk P :=\dim P-\dim \ker P-\dim \Indef P=\\=
  \dim\Dom P- \dim\Ker P
         =\dim\Im P-\dim\Indef P
         .\end{multline*}

{\sc Remark.}
Obviously, for any linear operators  $A:V\to W$, $B:W\to Y$,
\begin{align*}
&\graph(BA)=\graph (B)\,\graph (A);\\
&\Ker \graph (A)= \Ker A;\qquad \Im\graph (A)=\Im A\\
&\rk\graph (A)=\rk A.
\end{align*}

For $P:V\rr W$ we define the {\it pseudoinverse}
linear relation $P^\square:W\rr V$. It is the same subspace
$P\subset V\oplus W$ regarded as a subspace in $W\oplus V$.

For $P:V\rr W$ and $c\in \C^*$ we define the linear relation
$c\cdot P$ consisting of all vectors $v\oplus c w$,
where $v\oplus w$ ranges in $P$.

\smallskip

{\sc Remark.} Let $A$ be a  linear operator. Then
$c\cdot\graph(A)=\graph(c A)$.
For an invertible  linear operator $A:V\to V$, we have
$$(\graph A)^\square=\graph (A^{-1}).$$

{\bf 1.2. Category $\GA$.}
The objects of the category $\GA$ are finite-dimensional
linear spaces over $\C$. Set of morphisms
$\Mor(V,W)= \Mor_\GA(V,W)$ consists of all linear relations
$P:V\rr W$ and the formal element $\null_{V,W}$.

\smallskip

{\sc Remark.}
The dimension of $ P$ is an arbitrary number
$0,1,2,\dots, \dim V+\dim W$. The element $\null_{V,W}$ can not
be identified with
a linear relation.  \konets

 \smallskip

Let $P:V\rr W$, $Q :W\rr Y$ be linear relations.
If
\begin{align}
\Im P+ \Dom Q=W;\\
\Indef P\cap \Ker Q=0,
\end{align}
then the product $QP$ in the category $\GA$ coincides
with the product of linear relations. Otherwise,
$$QP=\null_{V,Y}.$$
The product of $\null$ and any morphism is $\null$.

\smallskip

{\sc Proposition 1.1.} (see \cite{Nerb}, 2.7)
 a) {\it For any linear spaces $V$, $W$, $Y$, $Z$
  and any morphisms
   $P\in \Mor(V,W)$, $Q\in\Mor (W,Y)$, $R\in \Mor(Y,Z)$,
   the associativity holds}
   $$(RQ)P=R(QP).$$

  b) {\it Let  $P\in \Mor(V,W)$,
  $Q\in\Mor (W,Y)$. If $QP\ne\null$,
  then}
$$\dim(QP)=\dim Q+\dim P -\dim W.$$

{\sc Remark.}  The group $\Aut(V)$ of automorphisms
of an object $V$ is the  group $\GL(V)$.     \konets

\smallskip

{\bf 1.3. Semigroup $\Gamma(V)$.}
Denote by $\Gamma(V)$ the subset in $\Mor(V,V)$
consisting of $\null_{V,V}$ and all linear relations
$R: V\rr V$ such that $\dim R=\dim V$.

\smallskip

By Proposition 1.1,  $\Gamma(V)$ is closed with respect to  the
multiplication.

Also, for a linear relation $R\in\Gamma(V)$
$$
\dim\Dom R+\dim\Indef R=\dim V;\qquad
\dim\Ker R +\dim \Im R=\dim V
$$

\smallskip

{\bf 1.4. Exterior algebras.} Let $V$ be a complex linear space.
Denote by $\Lambda(V)$ the exterior algebra of the space $V$.
Recall that $\Lambda(V)$
is the associative algebra with generators
$v$, where
$v$ ranges in $V$, and  the relations
\begin{align*}
& v\wedge w=-w\wedge v,\quad ;\\
&(\mu v_1+\nu v_2)\wedge w=
 \mu ({v_1}\wedge w) +  \nu ({v_2}\wedge w),
 \end{align*}
where $\mu,\nu\in \C$, $v,w\in V$ and
the sign $\wedge$ denotes the multiplication in $\Lambda(V)$.

We denote by $\Lambda^k V$
 the linear  subspace in $\Lambda(V)$
spanned by all vectors having the form
$$v_1\wedge v_2\wedge \dots \wedge v_k,\qquad v_j\in V.$$
The space $\Lambda^k V$ is called {\it the $k$-th exterior power}
of $V$.
If $e_1,\dots,e_n$ is a basis in $V$, then
the collection $e_{i_1}\wedge\dots \wedge e_{i_k}$,
where $i_1<i_2<\dots <i_k$, is a basis in $\Lambda^k V$.

\smallskip

   Let $A: V\to W$ be a linear operator.
We define the linear {\it operator of  change of the variables}
$$\lambda_{\change}(A):\Lambda V\to \Lambda W$$
 by
$$
 \lambda_{\change}(A) v_1\wedge \dots \wedge v_k=
                   Av_1\wedge \dots \wedge Av_k.
$$
If $A$ is an operator $V\to W$, and $B$ is an operator
$W\to Y$, then
$$\lambda_{\change}(B)\lambda_{\change}(A)
=\lambda_{\change}(BA).$$

The operators $\lambda_{\change}(A)$ preserve the degree $k$, and
hence we also obtain the operators
in the exterior powers
 $$\lambda^k_{\change}(A):\Lambda^k V\to \Lambda^k W.$$

\smallskip

{\bf 1.5. Fundamental representation of the category $\GA$.}
Let $S:V\rr W$ be a linear relation. Then there exist a basis
\begin{equation}
f_1,\dots,f_\alpha, g_1,\dots,g_\beta,
      h_1,\dots,h_\gamma
      \end{equation}
in $V$ and a basis
\begin{equation}
F_1,\dots,F_\mu, G_1,\dots, G_\beta,
      H_1,\dots, H_\nu
\end{equation}
in $W$ such that $S\subset V\oplus W$ is
spanned by the vectors
\begin{equation}
0\oplus F_1,\dots, 0\oplus F_\mu,\,\,\,
 g_1\oplus G_1, \dots,
 g_\beta\oplus G_\beta,\,\,\, h_1\oplus 0,
   \dots h_\gamma \oplus 0               .
   \end{equation}

{\sc Remark.} Thus,

-- the vectors $h_i$ form a basis in $\Ker S$;

-- the vectors $h_i$ and $g_j$ form a basis in $\Dom S$;

-- the vectors $F_k$ form a basis in $\Indef S$;

-- the vectors $F_k$, $G_j$ form a basis in $\Im S$.
 \hfill$\square$

 \smallskip

{\sc Remark.} Let $V=W$ and $S\in \Gamma(V)$.
 Then $\alpha=\mu$, $\gamma=\nu$.  \hfill$\square$

\smallskip

We define the linear operator
$$\lambda(S):\Lambda(V)\to\Lambda(W)$$
by
$$\lambda(S)f_1\wedge f_2\wedge \dots\wedge f_\alpha
     \wedge g_{i_1}\wedge\dots\wedge g_{i_k}
     =
      F_1\wedge F_2\wedge\dots\wedge F_\mu
     \wedge G_{i_1}\wedge\dots\wedge G_{i_k},
$$
and
$$\lambda(S)\xi=0$$
for all other basic vectors $\xi$ in $\Lambda V$.

\smallskip

{\sc Remark.} The bases (1.3), (1.4) are not uniquely determined by the
linear relation $S$. If we change the bases, then the operator
$\lambda(S)$  shall be multiplied by a nonzero constant.
\hfill$\square$

\smallskip

Let $A\in \GL(W)$, $B\in\GL(V)$. Then
$$\lambda(ASB)=\lambda_\change(A)
\lambda(S)\lambda_\change(B)
.$$
We also assume
$$\lambda(\null)=0.$$

{\sc Remark.} Let $S$ be a graph of a linear operator $A$.
Then $\lambda(S)=\lambda_\change(A)$. Nevertheless,
 $\lambda_\change(A)$ is a well defined operator in $\Lambda V$,
 the operator $\lambda(S)$ is defined up to a
 nonzero scalar factor.
By this reason, we preserve
the both notations  $\lambda(\cdot)$
and $\lambda_\change(\cdot)$, since their meanings
 slightly differ.\konets

\smallskip

{\sc Theorem 1.3.} (see \cite{Nerb}, II.7 or
 \cite{Ner5}, \S1) a) {\it Let $P:V\rr W$, $Q:W\rr Y$
be linear relations. Then
$$\lambda(Q)\lambda(P)=c(Q,P)\lambda(QP),$$
where $c(q,p)\in\C$.}

\smallskip

b) {\it
$c(q,p)\ne 0$ iff $QP\ne \null$.}

\smallskip

 {\sc Remark.} For the coordinateless
 definition of the operators $\lambda(P)$, see
 \cite{Ner2}, \cite{Nerb}

 \smallskip

{\sc Remark.} By the construction, the operator $\lambda(P)$ takes
homogeneous vectors to homogeneous vectors.
Thus we obtain the family of the operators
$$\lambda^k(P):\Lambda^k V\to\Lambda^{k+\dim P-\dim V} W.$$

\smallskip

{\bf 1.6. Fundamental representations of
the semigroup $\Gamma(V)$.}
The main tool below is the semigroup $\Gamma(V)$   defined in
1.3. Obviously, for $P\in\Gamma(V)$ we have
$$\lambda^m(P):\Lambda^m V \to  \Lambda^m V.$$
Thus we obtain the collection of the
projective representations
 $\lambda^m$ of the semigroup $\Gamma(V)$
 in the spaces $\Lambda^m V$.

 \smallskip

 {\sc Lemma 1.3.}  {\it Consider $P\in\Gamma(V)$.}

 \smallskip

 a) {\it $\lambda^m (P)\ne 0$
 iff $\dim \Indef P\le m \le \dim \Im P$}

 \smallskip

 b) {\it If $m=\dim\Indef P$, then $\rk \lambda^m(P)=1$,
 and
 $$\lambda^m(P)=\lambda^m(Q)$$
 for any $Q\in \Gamma(V)$ such that $\Indef Q=\Indef P$,
 $\Dom Q=\Dom P$.
In particular, we can choose  $Q=\Dom P\oplus \Indef P$.}

\smallskip

 c) {\it If $m=\dim\Im P$, then   $\rk \lambda^m(P)=1$,
 and
 $$\lambda^m(P)=\lambda^m(R)$$
 for any $R\in \Gamma(V)$ such that
 $\Im R=\Im P$, $\Ker R=\Ker P$. In particular, we can choose
 $R= \Ker P\oplus \Im P$}.

 \smallskip

 {\sc Proof.} Consider the canonical form described
 in 1.5. The conditions
 $\dim V=\dim W=\dim P$ imply
 $\alpha=\mu$, $\gamma=\nu$.
 Now all statements become obvious,
  see also \cite{Ner5}.    \konets

  {\sc Proposition 1.4.} {\it The map $\lambda(P)$ is a continuous
  map from the Grassmannian $\Gr_n$ of all $n$-dimensional
  subspaces in $V\oplus V$ to $\PMat(\Lambda V)$. In particular,
  its image is a closed subset
  in $\PMat(\Lambda V)$.}

 \bigskip

 {\large\bf 2. Hinges}

 \medskip
  \addtocounter{sec}{1}
\setcounter{equation}{0}

This Section contains the preliminaries on the hinges
with short explanations and sketches of proofs.
For more details, see \cite{Ner5}.

 In this Section, $V$
 is an $n$-dimensional complex linear space
 and $\Gr_n$ is the Grassmannian of all
 $n$-dimensional linear subspaces
 in $V\oplus V\simeq\C^n\oplus\C^n$.

 \medskip

 {\bf 2.1. Hinges.}
 A hinge
 $$\cP= (P_1,\dots,P_k):V\rrr V$$
 is a family of linear relations $P_j: V\rr V$
 such that $\dim P_j=\dim V=n$
 (hence $P_j\in \Gamma(V)$) and
 \begin{align}
 &\Ker P_j=\Dom P_{j+1},\qquad &\text{where}\quad j=1,2,\dots ,k-1,\\
 & \Im P_j=\Indef P_{j+1},\qquad & \text{where}\quad j=1,2,\dots ,k-1,\\
 & \Dom P_1=V,\\
 & \Im P_k=V,\\
 &
 P_j\ne \Ker P_j\oplus \Indef P_j
     ,\qquad & \text{where}\quad j=1,2,\dots ,k
 .\end{align}

 {\bf 2.2. Comments on the definition.}

 {\sc Remark.} The conditions (2.3)--(2.4)
 are the interpretation of the conditions
(2.1)--(2.2) for $j=0$ and $j=k$.
The condition (2.3) means that the first term $P_1$
of a hinge is  an operator.
The condition (2.4) means that
the last term $P_k$ is a linear relation pseudoinverse
to an operator.        \konets

\smallskip

{\sc Remark.} The graph of an {\it invertible} operator
is a hinge ($k=1$). The graph of a noninvertible operator
is not a hinge.                      \konets

\smallskip

{\sc Remark.} Let $A:V\to V$, $B:V\to V$
be linear operators such that
$$\Im A=\Ker B;\qquad \Ker A=\Im B.$$
Then
$$(\graph(A),\graph(B)^\square)$$
 is a hinge.
Any  hinge consisting of two terms
($k=2$) has this form.      \konets

\smallskip

{\sc Remark.} The condition (2.5)
is equivalent to the condition
$$\rk P_j>0.$$

\smallskip

{\sc Remark.} The condition (2.5) is technical.    For each hinge
$\cP$ we define the {\it completed hinge}
 \begin{equation}
 \widehat\cP:=(Q_0,P_1,Q_1,P_2,Q_2,\dots,P_k,Q_k)
 ,\end{equation}
 where
\begin{align}
&Q_0=V\oplus 0,\nonumber\\
& Q_j=\Ker P_j\oplus \Im P_j=\Dom P_{j+1}\oplus \Indef P_{j+1},\\
& Q_k=0\oplus V\nonumber
 ,\end{align}

Obviously, $\cP$ is uniquely determined by $\widehat \cP$.
Thus, the space of all hinges and the space of all completed hinges
coincide.\konets

\smallskip

{\sc Remark.}
We have
\begin{align*}
&V\supset\Ker P_1\supset \Ker P_2
            \supset\dots \supset \Ker P_k = 0,\\
&0\subset \Im P_1 \subset \Im P_2 \subset\dots\subset \Im P_k=V
.\end{align*}
By (2.5), $\Ker P_{j+1}\ne\Ker P_j$.
This implies $k\le n$.   \konets

\smallskip

{\sc Remark.}
We have
$$\Ker P_j\oplus \Indef P_j  \subset P_j
          \subset \Dom P_j\oplus \Im P_j.$$
The image of $P_j$ under the natural projection
 $$
  \Dom P_j\oplus \Im P_j \,
     \longrightarrow\, \bigl(\Dom P_j/\Ker P_j\bigr)
     \oplus\bigl( \Im P_j /\Indef P_j\bigr)
 $$
 is a graph of an invertible operator
 $$\Dom P_j\bigl/\Ker P_j \to \Im P_j \bigl/\Indef P_j .$$

 {\sc Remark.}
 Thus,   hinges can be defined in
 the following way. Consider
 two flags
 \begin{align*}
 0=Y_0\subset Y_1\subset \dots \subset Y_k = V,\\
 V=Z_0\supset Z_1\supset \dots \supset Z_k = 0 ,
 \end{align*}
 such that
 $$\dim Y_j/Y_{j-1}=\dim Z_{j-1}/ Z_j
              \qquad \text{for all}\quad j.$$
 For each $j$, we  fix an invertible linear operator
 $$A_j: Y_j/Y_{j-1}\to Z_{j-1}/ Z_j.$$
 By the previous remark, these data define a hinge.  \konets

 \smallskip

{\bf 2.3. Notation.}
 We denote by $\Hinge(V)=\Hinge_n$ the space of all hinges
 $\cP: V\rrr V$.

 We denote by $\Hinge^*(V)=\Hinge^*_n$ the space of all hinges
 $\cP=(P_1,\dots, P_k): V\rrr V$ defined up
 to the equivalence
 $$
 (P_1,\dots, P_k)  \sim (c_1 P_1,\dots,c_k P_k),
     \qquad \text{where} \quad c_j\in \C^*   .
 $$
 Considering 1-term hinges ($k=1$), we obtain
 $$\Hinge_n\supset\GL_n(\C); \qquad
	     \Hinge_n^*\supset\PGL_n(\C).$$

 \smallskip

 {\bf 2.4. Topology on $\Hinge_n^*$.}
 Below we define a structure of
 an irreducible smooth projective algebraic
 variety on $\Hinge_n^*$.
 In this Subsection, we define the topology on $\Hinge^*_n$.

\quad {\bf 2.4.a. Convergence of sequences $g_j\in \PGL_n$ to points
of $\Hinge_n^*$.}
 The sequence $g_j\in \PGL_n$ converges to
 $\cP=(P_1,\dots,P_k)\in\Hinge_n^*$, if there exist $k$ sequences
 \begin{equation}
 \beta^{(1)}_1, \beta^{(1)}_2,\beta^{(1)}_3,\dots;
 \qquad
 \beta^{(2)}_1, \beta^{(2)}_2,\beta^{(2}_3 \dots;
 \qquad
 {\bold \dots\dots};
  \qquad
  \beta^{(k)}_1, \beta^{(k)}_2,\beta^{(k)}_3,\dots
 ,\end{equation}
 ($\beta^{(\sigma)}_j\in\C^*$)
 such that   $\beta^{(\sigma)}_j g_j$
 converge to $P_\sigma$ in $\Gr_n$ for all $\sigma=1,\dots,k$.

 \smallskip

 {\sc Remark.} This implies
 $$\lim_{j\to \infty}
 \beta_j^{(\sigma)}/ \beta^{(\tau)}_j=\infty
 \qquad \text{for $\sigma>\tau$}.$$

 \smallskip

 {\sc Example.} Let
 $$g_j=\begin{pmatrix}
 4^j&0&0\\
 0& 2^j &0\\
 0&0& 2^{-j}
 \end{pmatrix}
 .$$
 Then we can choose
 $$\beta_j^{(1)}=4^{-j};\qquad
  \beta_j^{(2)}=2^{-j};\qquad
   \beta_j^{(3)}=2^{j}
 .$$
 For the sequences
 $$\mu_j^{(0)}=8^{-j};\quad
 \mu_j^{(1)}=2^{-3/2\, j};\quad
     \mu_j^{(2)}=1;
    \quad  \mu_j^{(3)}=4^{j}
  $$
   the limits $Q_j=\lim_{j\to\infty}\mu^{(\sigma)}_j g_j$ in $\Gr_n$
  also exist, but all the limits $Q_0$, \dots, $Q_3$ have rank 0.
  These limits are the elements $Q_j$ of the completed hinge, see
  (2.6).                     \konets

  \smallskip

 {\sc Lemma 2.1.} {\it
 Any sequence $g_j\in\PGL_n(\C)$
 contains a subsequence convergent in our sense.}

 \smallskip

 {\sc Proof.} We represent $g_j$ in the form
 $$g_j= A_j
 \begin{pmatrix}
                     u_j^{(1)}&0&\dots  \\
                      0& u_j^{(2)}&\dots   \\
		      \vdots&\vdots&\ddots
     \end{pmatrix}
    B_j,
    $$
    where $A_j$, $B_j$ are unitary matrices, and
    $$u_j^{(1)}\ge u_j^{(2)}\ge\dots >0$$
    Selecting a  subsequence, we can assume that

    1) the sequences $A_j$, $B_j$ are convergent

    2)$\exists \lim\limits_{j\to \infty}
       u_j^{(m)}/u^{(1)}_j=\alpha_m$ for all $m>1$.

 Obviously, $1\ge\alpha_2\ge\alpha_3\ge\dots\ge 0$.
 Consider $\tau$ such that $\alpha_{\tau-1}\ne 0$,
 $\alpha_\tau=0$. After the next selection of a subsequence,
 we can assume

  3) $\exists \lim\limits_{j\to \infty}
       u_j^{(m)}/u^{(\tau)}_j=\beta_m$
        for all $m>\tau$.
	Obviously,
$1\ge \beta_{\tau+1}\ge\beta_{\tau+2}\ge\dots\ge 0$.

 Then we repeat the same argument
 again, again, and again.

 Now we assume
 $\beta^{(1)}_j=\bigl(u_j^{(1)}\bigr)^{-1}$,
  $\beta^{(2)}_j=\bigl(u_j^{(\tau)}\bigr)^{-1}$,
  etc.
  \konets

  \smallskip
  Thus the space $\Hinge_n^*\setminus\PGL_n(\C)$
  is some kind of  boundary of $\PGL_n(\C)$.

  \smallskip

 {\bf 2.4.b. Convergence on the boundary.}
   A sequence
  $\cP^{j}=(P_1^j,\dots P_k^j)$ converges to
  a hinge $\cQ=(Q_1,\dots, Q_l)$, if for any
  $Q_u$, there exist $v=1,\dots,k$ and a sequence
  $\mu_j\in\C^*$ such that $\mu_j P^j_v$ converges
  to $Q_u$ in $\Gr_n$.

  \smallskip

 {\bf 2.4.c. Formal description of
the structure of a compact metric space on
   $\Hinge^*_n$.}
Consider the action
 of the group $\C^*$ on $\Gr_n$ given by
 $P\mapsto c\cdot P$.
 Fixed points of $\C^*$ in $\Gr_n$ are linear relations of
 rank 0.
    In other words, the fixed points have the form
    $$Q=\Ker Q\oplus \Indef Q.$$
 All other orbits of $\C^*$ in $\Gr_n$
 have trivial stabilizers.
 If $P\in \Gr_n$ is not a fixed point of $\C^*$-action,
 then the closure of the orbit
 $\C^*\cdot P$ consists of the orbit itself and
 of the pair of the points
 $$\Dom P\oplus \Indef P,\qquad
      \Ker P\oplus \Im P.$$
 For a hinge $\cP=(P_1,\dots,P_k)$,
 we define the  subset $\Omega(\cP)$
 in
 $\Gr_n$ by
 $$
 \Omega(\cP):=Q_0\,\cup\, \C^* P_1\,\cup\,
     Q_1\,\cup\, \C^* P_2\,\cup\, Q_2\,\cup\, \dots\,
         \, \cup\,  Q_{k-1}\,\cup\, \C^* P_k\, \cup\, Q_k
 ,$$
 where the points
 $Q_j=\Ker P_j\oplus\Im P_j=\Dom P_{j+1}\oplus\Indef P_{j+1}$
 are the elements of the completed hinge, see (2.6).

 \smallskip

  We emphasis that the closure
 of $\C^* P_j$ contains $Q_{j-1}$ and $Q_j$. Hence the subset
 $\Omega(\cP)$ is closed and connected.

 \smallskip

 Denote by $\Clos$ the space of all closed subsets
 of the Grassmannian $\Gr_n$.
 Consider an arbitrary metric $\rho$
 on $\Gr_n$ compatible with the topology.
 For $x\in\Gr_n$ and $A\in\Clos$, we define the
 distance
 $$\rho(x,A)=\min\limits_{y\in A} \rho(x,y).$$

 The {\it Hausdorff metric} in $\Clos$
 is defined by
 $$d(A,B)=\max\Bigl[
   \max\limits_{x\in A} d(x,B),\,\,
   \max\limits_{y\in B} d(y,A)\Bigr]
   .$$

{\sc Theorem 2.2} a)
{\it The image
of the embedding $\Hinge_n^*\to\Clos$
given by $\cP\mapsto \Omega(\cP)$ is
a compact subset in $\Clos$. Thus we obtain
a topology of a compact metric space on $\Hinge_n^*$.}

\smallskip

b) {\it  The group $\PGL_n(\C)$ is dense
in $\Hinge^*_n$.}

\smallskip

{\sc Remark.} The space of orbits of $\C^*$
on $\Gr_n$ is a nonseparated topological space.
Construction described above is the result of
an application
of the construction of the Hausdorff quotient described
in \cite{Ner4},\cite{Ner5}. In the algebraic geometry,
 there exist
also more delicate  constructions
of the Hilbert scheme quotient and
the Chow scheme quotient, see
\cite{BS},\cite{Kap}; in our case, these constructions
are equivalent to the Hausdorff quotient.             \konets

\smallskip

 {\bf 2.5. Orbits of the group
 $\GL_n(\C)\times \GL_n(\C)$
 on $\Hinge_n$.} The group $\GL_n(\C)\times\GL_n(\C)$
 acts on $V\oplus V\simeq \C^n\oplus\C^n$ in
 the obvious way.
 Hence it acts on the spaces $\Hinge_n$ and  $\Hinge^*_n$

 \smallskip

 {\sc Lemma 2.3.}
  {\it The group $\GL_n(\C)\times \GL_n(\C)$
 has $2^n$ orbits on the space $\Hinge_n$.
  These orbits are
 enumerated by the
  number $k=1,2,\dots, n$ and the  positive
 numbers}
    $$\alpha_1=\rk P_1,\dots,  \alpha_k=\rk P_k,
    \qquad\text{where}\quad \alpha_1+\dots+\alpha_k=n.$$

 \smallskip

 {\sc Proof.} Obvious.         \konets

 \smallskip

We
 denote these orbits by
$$\cO[\alpha]=\cO[\alpha_1,\dots,\alpha_k].$$

Fix a basis $e_1,\dots, e_n\in V$.
For a given collection  $\alpha_1,\dots,\alpha_k$,
 we define
the {\it canonical hinge}
$\cP_{\alpha_1,\dots,\alpha_k}\in \cO[\alpha_1,\dots,\alpha_k]$
by
$$\cP_{\alpha_1,\dots,\alpha_k}=(P_1,\dots,P_k):V\rrr V,$$
where the linear relation $P_j\subset V\oplus V$ is
spanned by the vectors
\begin{align*}
&0\oplus e_\sigma,
           \qquad &\text{where}\quad
           \sigma\le\alpha_1+\dots + \alpha_{j-1},\\
& e_\tau\oplus e_\tau,
           \qquad &\text{where}\quad
            \alpha_1+\dots + \alpha_{j-1}<
                      \tau             \le
             \alpha_1+\dots + \alpha_{j} ,\\
&   e_\mu\oplus 0,   \qquad &\text{where}\quad
              \mu>\alpha_1+\dots + \alpha_{j}.
              \end{align*}

{\sc Remark.} Denote by $\cO^*[\alpha]$ the image of
$\cO[\alpha]$ in $\Hinge_n^*$. We have
 $$\dim\cO[\alpha]=n^2 ;
    \qquad \dim \cO^*[\alpha_1,\dots,\alpha_\tau]=n^2-\tau$$

 {\bf 2.6. Alternative.}

{\sc Theorem 2.4.} {\it Let $\cP=(P_1,\dots,P_k):V\rrr V$
 be a hinge.
Fix $m=0,1,\dots,n$.
Consider the family of operators
\begin{equation}
\lambda^m(P_1),\lambda^m(P_2),\dots, \lambda^m(P_k):
\Lambda^m V\to \Lambda^m V
.\end{equation}

Then there are only two possibilities.}

\smallskip

 1) {\it There exists a unique $j$ such that $\lambda^m(P_j)\ne 0$.}

\smallskip

  2) {\it There exists $j$ such that $\lambda^m(P_j)\ne 0$,
 $\lambda^m(P_{j+1})\ne 0$, and $\lambda^m(P_\tau)= 0$
 for all $\tau\ne j,j+1$. In this case,
 $\lambda^m(P_j)$ and
 $\lambda^m(P_{j+1})$ have rank 1 and coincide up to a nonzero factor.
 They also coincide with $\lambda(Q_j)$, where
 $Q_j=\Ker P_j\oplus\Im P_j$ is the term of the
 completed hinge $\widehat \cP$, see {\rm (2.6)}.}

 \smallskip

 {\sc Proof.} This is a  consequence of Lemma 1.4. \konets

 \smallskip

 Now for each hinge $\cP=(P_1,\dots,P_k):V\rrr V$
 and for each $m=0,1,\dots,n$,
 we  define the operator
 $$\lambda^m(\cP):\Lambda^m V\to \Lambda^m V$$
 as the unique nonzero term of the sequence (2.9).
 By the construction, this operator
 is  defined up to a nonzero factor.
   \smallskip

   {\sc Remark. } Obviously, for any $g_1,g_2\in\GL_n(\C)$,
   $$\lambda^m(g_1\cP g_2)=
      \lambda^m_{\change}(g_1)
         \lambda^m(\cP)\lambda^m_{\change}(g_2).$$

 \smallskip

 {\bf 2.7. Example:
 the operators $\lambda^m(\cP)$ for canonical hinges.}
   Consider the canonical hinge
 $\cP_\alpha=\cP_{\alpha_1,\dots, \alpha_k}$.
  We intend to describe the operator
$\cL^m:=\lambda^m(\cP_{\alpha})$.
  Assume
 $$u_j= \alpha_1+\dots+\alpha_{j-1}$$

Let $u_j\le m\le u_{j+1}$,
let $s:=m-u_{j}$ and $u_j<l_1< \dots <l_s\le u_{j+1}$.
Then for each element
$$h
  =e_1\wedge e_2\wedge e_3\wedge\dots\wedge e_{u_j}\wedge
        e_{l_1}\wedge \dots\wedge e_{l_s}$$
        of the standard basis,
        we have
$$\lambda^m(\cP_\alpha) h=h;$$
and $\lambda^m(\cP_\alpha)$
 annihilates all other elements
$e_{t_1}\wedge\dots\wedge e_{t_m}$ of the standard basis
in $\Lambda^m V$.

 \smallskip

 {\bf 2.8. The projective embedding of $\Hinge^*_n$.}
In 2.6,
 for any $\cP\in\Hinge_n$, we constructed the family
of nonzero linear operators
\begin{equation}
\lambda^\circ(\cP):=
 \bigl( \lambda^1(\cP),   \lambda^2(\cP),\dots,
     \lambda^{n-1}(\cP)
      \bigr)
      \end{equation}
defined up to nonzero factors.
Consider two hinges
$$\cP=(P_1,\dots, P_k);\qquad \cR= (c_1 P_1,\dots,c_k P_k),
 \quad \text{where}\quad c_j\in \C^*.$$
Obviously,
the operators
 $\lambda^m(\cP)$ and $ \lambda^m(\cR)$ coincide up to
a nonzero factor.

 Thus we obtain the map
 $$
\lambda^\circ:\Hinge^*_n\to
     \text{\huge $\times$}_{m=1}^{n-1}
     \PMat(\Lambda^m V)              .
 $$

Consider also the map
$$\lambda^\circ:\PGL_n(\C)
  \to \text{\huge $\times$}_{m=1}^{n-1}\PMat(\Lambda^m V)$$
given by
\begin{equation}
\lambda^\circ(g):=
   \bigl( \lambda^1_{\change}(g),
      \dots,\lambda^{n-1}_{\change}(g)\bigr)
      .\end{equation}

{\sc Theorem 2.5.} (\cite{Ner3},\cite{Ner5})
{\it The map {\rm(2.10)}  from $\Hinge_n^*$ to
$ \text{\huge $\times$}_{m=1}^{n-1}\PMat(\Lambda^m V)$
         is continuous.}

\smallskip

{\sc Corollary 2.6.}
a) {\it The image
of the map
{\rm (2.10)}
is compact.}

\smallskip

b) {\it The image of the map {\rm(2.10)}
is the closure of the $\lambda^\circ\bigl(\PGL_n(\C)\bigr)$.}

\smallskip

c) {\it The space $\Hinge_n^*$
is an irreducible projective variety.}

\smallskip

{\sc Sketch of proof of Theorem 2.5.}
We shall prove the implication
 \begin{equation}\text{${g_j}\in\GL_n(\C)$ converges to $\cP$}
\,\, \Longrightarrow \,\,
 \text{$\lambda^\circ(g_j)$ converges
to $\lambda^\circ(\cP)$}
.\end{equation}

Let us represent $g_j$ in the form
$$g_j= A_j D_j B_j,$$
where $A_j$, $B_j$ are unitary matrices
and $D_j$ are  diagonal matrices with
the decreasing eigenvalues. We can assume that
the sequences $A_1, A_2,\dots$ and
$B_1,B_2,\dots$ are convergent. Thus
the question is reduced to the case
$A_j=1$, $B_j=1$. For this case,
the statement can easily  be checked.

\smallskip

 Corollary 2.6 implies the following consequence.

 \smallskip

{\sc Corollary 2.7.}
 a) {\it The space $\Hinge_n^*$ is an irreducible smooth
 projective variety.}

\smallskip

  b) {\it The variety  $\Hinge_n^*$
coincides with the Semple complete collineation space.}

\smallskip

c) {\it  The variety  $\Hinge_n^*$ coincides with
 the De Concini--Procesi {\rm \cite{DCP}}
  compactification of $\PGL_n(\C)$.}

 \smallskip

By the definition (see \cite{Sem2}),
 the Semple complete collineation space
is the closure of the image of the map (2.11).
By the Semple theorem, the complete collineation
space is a smooth projective variety, and by \cite{DCP}
it coincides with
 the De Concini--Procesi compactification,
 see also below 2.13.

\smallskip

{\bf 2.9. Semigroup of hinges.}
Let us define one more variation of the space $\Hinge_n$.
Denote by $\widetilde\Hinge_n$ the set of all elements
of $\text{\huge$\times$}_{m=0}^n\, \Mat(\Lambda^m V)$
having the form
$$\cA=(A_0,\dots,A_n)=
\bigl(c_0\cdot\lambda^0(\cP), c_1\cdot\lambda^1(\cP),
   c_2\cdot\lambda^2(\cP),\dots,c_n\cdot\lambda^n(\cP)\bigr),$$
where $\cP$ is a hinge and $c_0,\dots,c_{n}\in\C$.
We say that $\cA$ {\it lies
over $\cP$.} If all $c_j$ are nonzero,
we say that $\cA$ is {\it nondegenerated}.

\smallskip

{\sc Proposition 2.8.} {\it $\widetilde\Hinge_n$
 is a subsemigroup
in $\text{\huge$\times$}_{m=0}^n \Mat(\Lambda^m V)$.}

\smallskip

This easily follows from Theorem 2.5.

Now we intend to give
 a constructive description of the product in
$\widetilde\Hinge_n$.

Let $\cR=(R_1,\dots,R_s)$ be a family of
linear relations $V\rr V$, and $\dim R_j=n$. We say that
$\cR$ is a {\it weak hinge} if for each $j$
\begin{align*}
& \Ker R_j\supset\Dom R_{j+1},\\
& \Im R_j  \subset \Indef R_{j+1}.
\end{align*}

{\sc Remark.}
 Let $\widehat\cP=(Q_0,P_1,Q_1,P_2,\dots,P_k,Q_k)$ be a
completed hinge (see (2.6)). Then
any subcollection of $\widehat\cP$   is a weak hinge,
and  each weak hinge can be obtained in this way.      \konets

\smallskip

For any weak hinge $\cR=(R_1,\dots,R_s)$ and any
$m=0,1,\dots,n$, we intend to construct
the canonical operator
$$\lambda^m(\cR):\Lambda^m V\to\Lambda^m V$$
defined up to a scalar factor.
For this, we consider the sequence
$$\lambda^m(R_1),\lambda^m(R_2),\dots, \lambda^m(R_s).$$
If this family contains a nonzero term
$\lambda^m(R_j)$, then
$$\lambda^m(\cR):= \lambda^m(R_j).$$
Otherwise,
 $$ \lambda^m(\cR):=0.$$

 {\sc Theorem 2.9.}
 (\cite{Ner3}, \cite{Ner5})
  {\it Let $\cR=(R_1,\dots,R_s)$,
 $\cT=(T_1,\dots,T_\tau)$ be weak hinges.
 Then the  family of all
  $$T_iR_j\ne\null$$
  is a weak hinge.}

  \smallskip

  {\sc Sketch of proof.}    Let $T_iR_j\ne\null$.
  It can  easily be checked that the segments
  \begin{align}
  \dim \Indef R_j\le m\le \dim\Im R_j,\\
  \dim \Indef T_i\le m\le \dim\Im T_i
   \end{align}
   have nonzero intersection. After this remark,
    Theorem 2.8 can  easily be checked.    \konets

  \smallskip

  {\sc Theorem 2.10.}
   (\cite{Ner3}, \cite{Ner5}){\it For each $m$,
  $$\lambda^m(\cT)\lambda^m(\cR)= c\cdot\lambda^m(\cT\cR),$$
  where $c\in\C^*$.}

  \smallskip

  {\sc Proof.} Assume  $m$ satisfies the
   equations  (2.13)--(2.14). Then
  $$ \lambda^m(\cT)=\lambda^m(T_i);
     \qquad \lambda^m(\cR)=\lambda^m(R_j).$$
     Thus,
     $$\lambda^m(\cT)\lambda^m(\cR)=
     \lambda^m(T_i)\lambda^m(R_j)=
     \lambda^m(T_iR_j)=\lambda^m(\cT\cR). \qquad \square $$

     \smallskip

     Thus, for $\cA$ lying over a weak hinge $\cP$ and $\cB$ lying over
     a weak hinge $\cQ$, the product $\cA\cB$ lies over $\cP\cQ$.

     \smallskip

  {\bf 2.10. Canonical embedding
   $\Hinge_n\to \widetilde\Hinge_n$.}
  Let $\cP\in\Hinge_n$. The construction  2.6 defines
  the operators $\lambda^m(\cP)$ up to  nonzero factors.
  We intend to define these operators in a canonical way.

  Fix $\alpha=(\alpha_1,\dots\alpha_k)$.
  Consider $\cP=(P_1,\dots, P_k)\in\cO[\alpha]$.
  In particular, $\dim \Im P_j=\alpha_1+\dots+\alpha_j$,
  and $\dim \Indef P_j=\alpha_1+\dots+\alpha_{j-1}$

  We have the family of the operators
  $$\lambda(P_1),\,\lambda(P_2),\,\dots,
               \lambda(P_k):\Lambda(V)\to\Lambda(V)$$
               defined up to  nonzero factors.
  We have
  $$\lambda^m(P_j)\ne 0
     \qquad \text{iff}\quad \dim\Indef P_j\le m\le\dim\Im P_j,$$
  and
  \begin{align*}
  \lambda^{\alpha_1}(P_2)&=
          c_1\cdot \lambda^{\alpha_1}(P_1),\\
  \lambda^{\alpha_1+\alpha_2}(P_3)&
        =c_2\cdot \lambda^{\alpha_1+\alpha_2}(P_2),\\
   \lambda^{\alpha_1+\alpha_2+\alpha_3}(P_4)&
        =c_3\cdot \lambda^{\alpha_1+\alpha_2+\alpha_3}(P_3),
\end{align*}
etc.
  The linear relation $P_1$ is
  a graph of some linear operator $A$.
  Thus the operator
  $\lambda(P_1):=\lambda_\change(A)$ is well defined.
  After this we define the operator $\lambda(P_2)$ by the condition
  $c_1=1$,  then we define the operator
$\lambda(P_3)$ by the condition $c_2=1$  etc.

We denote by $\cL^m(\cP):\Lambda^mV\to\Lambda^m V$
the unique nonzero operator among
$$\lambda^m(P_1),\dots,\lambda^m(P_k).$$
We denote by $\cL(\cP)$ the collection $(\cL^0(\cP),\cL^1(\cP),\dots,\cL^n(\cP)$.
Thus we obtain the embedding
$$\cL: \Hinge_n\to
  \widetilde\Hinge_n.  $$

{\sc Remark.} For the canonical hinge
$\cP_\alpha$,  the family of the operators
 $\cL^m=\cL^m(\cP_\alpha)$
 in the exterior powers $\Lambda^m(V)$ was described in
 2.7.                   \konets

\smallskip

 {\sc Lemma 2.11.} {\it For any $g_1,g_2\in\GL_n(\C)$,}
   $$\cL^m(g_1\cP g_2)=
      \lambda_{\change}^m(g_1)
         \cL^m(\cP)\lambda_{\change}^m(g_2).$$

\smallskip

{\sc Proof.} Indeed, the multiplication by
 $\lambda_\change(g)$
 does not change the "gluing conditions"
 $c_1=c_2=\dots=1$.      \konets

 \smallskip

 {\bf 2.11. Reduced hinge semigroup.}
Consider the image $\cL(\Hinge_n)$ of the embedding $\cL$ described
in the previous subsection.     Denote by $\overline\Hinge_n$
the closure of this image.

The set    $\overline\Hinge_n$
is the union of $2^{n-1}$ affine algebraic varieties
$\overline{\cL(\cO_\alpha)}$, the dimension of all these varieties is
$n^2$.

The set $\overline\Hinge_n$ admits the following explicit description.
Let $\cR=(R_1,\dots,R_k)$ be a weak hinge.
Consider the family of  operators
$$c_j\cdot \lambda(R_j):\Lambda V\to \Lambda V,$$
where $c_j\in \C^*$. We say that the family is {\it well glued}
if

1. the condition $\lambda^m(R_j)\ne 0$, $\lambda^m(R_{j+1})\ne 0$
(or, equivalently, $\Ker(R_j)=\Dom(R_{j+1})$, and
their dimension is $m$)
implies
      $$c_j\cdot\lambda^m(R_j)=c_{j+1}\cdot\lambda^m(R_{j+1});$$

2.
if $R_1$ is an operator, then $c_1\lambda(R_1)=\lambda_\change(R_1)$.

\smallskip

The set $\overline\Hinge_n$
 coincides with the set of all well-glued families.

Obviously, the multiplication
of hinges preserves the glueing condition.
This implies the following statement

\smallskip

{\sc Proposition 2.12.} {\it The set $\overline\Hinge_n$ is a subsemigroup
in the semigroup $\widetilde\Hinge$.}

\smallskip

{\bf 2.12. Representations of the semigroup of hinges.}
Recall the construction of irreducible
polynomial finite-dimensional
representations of $\GL_n(\C)$.

Consider a collection of integers
\begin{equation}
\nu_1\ge\nu_2\ge\dots\ge\nu_n\ge 0.
\end{equation}
We  call such collections  {\it signatures}.
Denote by $\pi$ the standard representation of $\GL_n(\C)$
in $V=\C^n$. Fix a basis $e_1,\dots,e_n$  in $\C^n$.
Denote by $\xi_s$ the vector
$$\xi_s=e_1\wedge_2 \wedge\dots \wedge e_s\in \Lambda^s V.$$
Consider the space
$${\frak H}_\nu:=
  \bigotimes\limits_{j=1}^{n}
     \Lambda^j(V)^{\otimes (\nu_j-\nu_{j+1})}
$$
(we assume $\nu_{n+1}=0$). Consider the representation
${\frak r}_\nu$
of the group $\GL_n(\C)$ in the space ${\frak H}_\nu$
given by
$${\frak r}_\nu(g)=
  \bigotimes\limits_{j=1}^{n}
     \lambda_{\change}^j(g)^{\otimes (\nu_j-\nu_{j+1})} .
$$

Consider also the vector
$$\Xi_\nu:= \bigotimes\limits_{j=1}^{n}
     \xi_j^{\otimes (\nu_j-\nu_{j+1})}\in  {\frak H}_\nu.$$
Denote by $H_\nu$ the $\GL_n(\C)$-cyclic span
of the vector $\Xi_\nu$. We denote by $\rho_\nu$ the
representation of the group $\GL_n(\C)$ in the space
$H_\nu$.

It is well known (see, for instance \cite{Zhe}),
that the representations $\rho_\nu$ are irreducible
and all the polynomial irreducible representations
of the group $\GL_n(\C)$ have this form.

Let us define the representations $\rho_\nu$
of the semigroup $\widetilde\Hinge_n$.
Let
$${\cA}:=(A_0,A_1,A_2,\dots,A_n)\in \widetilde\Hinge_n; \qquad
   A_j\in \Mat(\Lambda^j V).$$
We define the operator
${\frak r}_\nu(\cA)$ in ${\frak H}_\nu$
by
$$  {\frak r}_\nu(\cA):
=
   \bigotimes\limits_{j=1}^{n}
     A_j^{\otimes (\nu_j-\nu_{j+1})}.
$$

{\sc Lemma 2.13.} {\it The subspace $H_\nu\subset {\frak H}_\nu$
is invariant with respect to the operators}
${\frak r}_\nu(\cA)$.

\smallskip

{\sc Proof.} Assume ${\frak r}_\nu(\cA)\ne 0$. Then the operator  ${\frak r}_\nu(\cA)$
depends (up to a scalar factor) only on the
hinge $\cP$ lying under $\cA$. But any hinge can be approximated by
elements of $\GL_n(\C)$.    \konets

\smallskip

We define the operator
$$\rho_\nu(\cA)$$
as the restriction of the operator   ${\frak r}_\nu(\cA)$
to the subspace $H_\nu$.
Obviously, $\rho_\nu$ is a linear
 representation of the semigroup
$ \widetilde\Hinge_n$:
  $$\rho_\nu(\cA)\rho_\nu(\cB)=\rho_\nu(\cA\cB)$$

{\sc Lemma 2.14.} {\it Let $\cA$ be a nondegenerated
element of $ \widetilde\Hinge_n$ lying
over the canonical hinge
 $\cP_{\alpha_1,\dots,\alpha_k}$. Then
$$\rho_\nu(\cA)\Xi_\nu=c\cdot \Xi_\nu,$$
where $c\in \C^*$.}

\smallskip

{\sc Proof.} Indeed,
$$\lambda^m(\cP_{\alpha_1,\dots,\alpha_k})
e_1\wedge\dots\wedge e_m=e_1\wedge\dots\wedge e_m,$$
see 2.7. This implies the required statement.    \konets

\smallskip

{\sc Corollary 2.15.} {\it Let $\cA$
be a nondegenerated element of $\widetilde\Hinge_n$.
 Then $\rho_\nu(\cA)\ne 0$.}

 \smallskip

{\sc Proof.} Each hinge  $\cQ$ is of the form $g_1\cP_\alpha g_2$,
there $g_1$, $g_2\in \GL(V)$.
Thus,
$$\rho_\nu(\cA)=
   c\cdot  \rho_\nu(g_1)\rho_\nu(\cP_\alpha)\rho_\nu(g_2).$$
The first and the third factors are invertible
and the middle factor is nonzero.    \konets

\smallskip

{\bf 2.13.  Extension of reducible representations.}
Let $\zeta$ be a finite-dimensional
polynomial
 representation
of $\GL_n(\C)$ in the space $Z$.
Then
$$\zeta=\oplus \rho_{\nu^{(j)}};\qquad Z=\oplus H_{\nu^{(j)}},$$
where $\nu^{(1)},\nu^{(2)}\dots $  are collections
 of signatures
satisfying (2.15), and $H_\nu$ are the corresponding spaces.
We define the representation
$\zeta(\cA)$ of the semigroup  $ \widetilde\Hinge_n$
by
$$\zeta(\cA)=\oplus  \rho_{\nu^{(j)}}(\cA)$$

{\bf 2.14. Identification of $\Hinge^*_n$ with
the De Concini--Procesi compactification
of $\GL_n(\C)$.}
Consider an {\it irreducible} representation $\rho_\nu$
of the group $\GL_n(\C)$.
Let $\cA\in\widetilde\Hinge_n $ lies over a hinge $\cP$.
Then the operator $\rho_\nu(\cA)$ is determined up to a factor by
the underlying hinge $\cP$. Moreover, it is determined
by $\cP$ considered as an element of $\Hinge^*_n$.
By Corollary 2.15, $\rho_\nu$ determines the map
$$\Hinge_n^*\to \PMat (H_\nu)$$
It can easily be checked that this  map  is continuous.

\smallskip

{\sc Theorem 2.16.} a) {\it Let $\rho_\nu$ be an irreducible
polynomial representation of $\GL_n(\C)$.
Then the image of $\Hinge^*_n$ in $\PMat(H_\nu)$
coincides with the projective compactification
$[\GL_n]_{\rho_\nu}$ defined in}  0.1.

\smallskip

b) {\it If $\nu_1>\nu_2>\dots>\nu_n$, then
the map $\rho_\nu:\Hinge_n^*\to\PMat(H_\nu)$
is an embedding.}

\smallskip

The statement b) identifies $\Hinge_n^*$ with
the De Concini--Procesi construction,
\cite{DCP}.

\bigskip

{\large \bf  3. Sea Urchin.}

\nopagebreak

\medskip
  \addtocounter{sec}{1}
\setcounter{equation}{0}

In this Section, we show that
the hinge language is sufficient
for a description of all the projective compactifications of $\GL_n(\C)$.

\smallskip

{\bf 3.1. Meromorphic matrices.}
Denote by $D_\epsilon$ the disk $|z|<\epsilon$
on $\C$.
Denote by $\dot D_\epsilon$ the punctured disk
$0<|z|<\epsilon$
on $\C$.

Denote by $\GL_N(\cF)$
the group of germs of holomorphic maps
$D_\epsilon\to \GL_N(\C)$.
 The elements of this group are
$N\times N$ matrices

$$\gamma(z)=\begin{pmatrix}
  \gamma_{11}(z)&\dots & \gamma_{1N}(z)\\
   \vdots&\ddots&\vdots  \\
  \gamma_{N1}(z)&\dots & \gamma_{NN}(z)
  \end{pmatrix}    ,
$$
where
$\gamma_{ij}(z)$ are functions holomorphic in a neighborhood
of $0$, and $\gamma(0)$ is invertible.

 We say that a map
 $$\gamma:\dot D_\epsilon \to \GL_N(\C)$$
 is a {\it meromorphic family} if all matrix elements
 $\gamma_{ij}$ are
 holomorphic functions in some punctured
 disk $\dot D_\epsilon$
 with poles or removable singularities
 at 0, and $\gamma(z)$ is invertible
 for $z$ lying in some punctured disk $\dot D_{\epsilon'}$.

 We define the order $\ord(\gamma)$
  of the pole of $\gamma$ as the maximal
 order of poles at 0 of the matrix elements $\gamma_{ij}$.
 In this definition, we admit a negative order
 of a pole (a function has a pole of
 negative order $-k$
 at a point 0, if it has the zero of order $k$ at $0$).

 The value $\ord(\gamma)$ coincides with the minimal $k$
 such that the map
 $$z\mapsto z^k \gamma(z)$$
 from $D_\epsilon$ to $\Mat(\C^N)$ is holomorphic.

 \smallskip

 {\bf 3.2. Exponents  of meromorphic families.}

 {\sc Lemma 3.1.} a) {\it Any meromorphic family
 $\gamma(z):\dot D_\epsilon\to \GL_n(\C)$
 can be represented in the form
 \begin{equation}
 \gamma(z)=a(z)
          \begin{pmatrix} z^{-m_1}&0&0&\dots\\
                          0&z^{-m_2}&0&\dots\\
                          0&0&z^{-m_3}&\dots\\
                          \vdots&\vdots&\vdots&\ddots
          \end{pmatrix}
          b(z),
 \end{equation}
 where $a(z), b(z)\in\GL_n(\cF)$ and}
 \begin{equation}
 m_1\ge m_2\ge\dots \ge m_n
\end{equation}

 b) {\it Consider the $j$-th exterior power
 $\lambda^j_{\change}(\gamma(z))$
  of the matrix $\gamma(z)$.
 Then
 $$\ord \lambda^j_{\change} (\gamma(z))
  =m_1+m_2+\dots+m_j$$
 In particular, the numbers $m_j$ are uniquely
 determined by the function $\gamma(z)$.}

 \smallskip

 {\sc Proof. } a) It is sufficient to apply
 the Gauss elimination algorithm.

 b) Obvious.                \konets

 \smallskip

 We say that the numbers $m_i$ are the {\it exponents}
 of the meromorphic family $\gamma(z)$.
 We fix the notation
 $${\bold m}:=(m_1,\dots, m_n)$$
  for the exponents. We also
 define the numbers
 $$k_1>k_2>\dots >k_\tau$$
to be all the pairwise different exponents $m_i$.
 We denote by $\alpha_j$ the number of copies
 $k_j$ in the collection $(m_1, m_2,\dots,m_n)$.
 We say that $\alpha_j$, $k_j$ are the {\it numbers associated
 with} $m_i$.

 \smallskip

{\it Starting this place, the sense of numbers
 $m_i$, $\alpha_j$, $k_j$ is fixed.}

 \smallskip

 {\bf 3.3. Limits of meromorphic families
  in the $\Hinge_n$.}
 Let $\gamma(z)$ be a meromorphic family
 of $n\times n$ matrices. Let $m_1,\dots,m_n$
 be its exponents, and let $k_1,\dots, k_\tau$,
 $\alpha_1,\dots,\alpha_\tau$ be the associated numbers.

 \smallskip

 {\sc Proposition 3.2.} a) {\it For each $j$
  the following limit in $\Gr_n$
 \begin{equation}
 P_j:=\lim\limits_{z\to 0}
   z^{k_j} \graph(\gamma(z))
 \end{equation}
 exists, and   $\rk P_j=\alpha_j$.}

   b) {\it The collection
   $\cP^\gamma:=(P_1,\dots,P_\tau)$
   is a hinge.
   }

   \smallskip

   {\sc Remark.} In 2.4, we defined the convergence
   of sequences $g_j\in\PGL_n(\C)$
   to elements of $\Hinge_n^*$. It is impossible
   to define
   a convergence of sequences $g_j\in\GL_n(\C)$ to  points of $\Hinge_n$.
   But in Proposition 3.2, we
    consider holomorphic families instead of sequences.

   \smallskip

 {\sc Proof.} Let
 $\gamma(z)$ be the diagonal matrix with the eigenvalues
 $z^{-m_1},\dots,z^{-m_n}$.
  In this case, the limit of $\gamma(z)$ is the canonical
  hinge $\cP_\alpha$. For a general family (3.1),
  the limit is
  $$ a(0)\cP_\alpha b(0)$$

 \smallskip



 {\bf 3.4. Limits of meromorphic families in
  $\overline\Hinge_n$.}
 Let $\gamma$ be a meromorphic family,
  let $m_j$ be its exponents.

 \smallskip

 {\sc Proposition 3.3} a) {\it For each $j=0,1,\dots, n$
 there exists the nonzero limit
 \begin{equation}
 \cL^j=\lim\limits_{z\to 0}
     z^{m_1+\dots+m_j} \lambda^j_{\change}(\gamma(z))
     \end{equation}
 in $\Mat(\Lambda^j V)$.}

\smallskip

 b) {\it The collection $\cL_j$ coincides with   the
    the collection $\cL^j(\cP^\gamma)$,
    where $\cP^\gamma$ is the hinge constructed
     in Proposition {\rm 3.2}
    and the operators $\cL^j(\cP)$
    were described in} 2.10.

\smallskip

   {\sc Example.} Consider
   the diagonal matrix
   $\delta_{\bold m}(z)$ with the eigenvalues
   $z^{-m_1},\dots,z^{-m_n}$.
   Then $\cP^{\delta_{\bold m}}$
   is the canonical hinge $\cP_{\alpha}$.
   The operator
    $\lambda^j_{\change}(\delta_{\bold m}(z))$
    is the diagonal operator in the basis
    $e_{p_1}\wedge\dots\wedge e_{p_j}$.
    The eigenvalues are $z^{-m_{p_1}-\dots-m_{p_j}}$.
    The maximal absolute value of the eigenvalues
    is $|z|^{-m_1-\dots-m_j}$.

    The eigenspace $W_j\subset \Lambda^j V$
     corresponding
    to the maximal eigenvalue is spanned by
    the vectors
     \begin{equation}
     e_1\wedge e_2\wedge \dots\wedge e_u\wedge
    e_{q_1}\wedge\dots\wedge e_{q_{j-u}}       ,
    \end{equation}
    where $u$ is the largest $i$ such that
    $m_i> m_j$ and $m_{q_1}=\dots=m_{q_{j-u}}=m_j$.
    Obviously,
    $$\cL^j=\lim\limits_{z\to 0} z^{m_1+\dots+m_j}
     \lambda_{\change}(\delta_{\bold m}(z))$$
     is the identical operator on the subspace spanned
     by the vectors (3.5), and $\cL^j$ annihilate
     all other basic vectors.
     Thus, $\cL^j$ coincides with the operator
     $\cL^j(\cP_\alpha)$ described in 2.7.  \konets

     \smallskip

     {\sc Proof of Proposition 3.3.}
     Let us represent $\gamma(z)$ in the form (3.1).
     Then, by the example given above,
     $$\cP^\gamma=a(0)\cP_\alpha b(0);\qquad
       \cL^j=\lambda^j_{\change}\bigl(a(0)\bigr)
           \cL^j(\cP_\alpha)
	   \lambda^j_{\change}\bigl(b(0)\bigr).$$
     Now we apply Lemma 2.11.    \konets


\smallskip

 {\bf 3.5. Limits of meromorphic families in irreducible
  representations.}
  Let $\rho_\nu$ be an irreducible holomorphic representation
  of $\GL_n(\C)$. We extend $\rho_\nu$
  to the representation of
  $\widetilde\Hinge_n$  by the procedure 2.12.

  {\sc Lemma 3.4.}  {\it Let $\gamma$
   be a meromorphic family
  of $n\times n$ matrices,
  let $m_i$ be its exponents.
  Let $\cL :\Hinge_n\to\widetilde\Hinge_n$
     be the embedding defined in} 2.10.

  \smallskip

    a) {\it There exists the nonzero limit}
    \begin{equation}
  [\rho_\nu(\gamma)]
  :=  \lim\limits_{z\to 0} z^{\sum m_i\nu_i} \rho_\nu(\gamma(z))
    .\end{equation}

   \smallskip

    b) {\it This limit coincides with the operator}
    $$\rho_\nu
       \bigl(\cL(\cP^\gamma)\bigr).$$

 {\sc Proof.} Obvious.    \konets

   \smallskip

 {\bf 3.6. Change of the variable $z$
  in meromorphic families.}
 We intend to consider a limit of a
 meromorphic curve independently
 on its parametrization.

 \smallskip

 a) If we change the variable $z$ by the formula
  $$z\mapsto z+c_2z^2+c_3z^3+\dots,$$
  then nothing  changes. The exponents $m_j$,
   the hinge $\cP^\gamma$, the limit (3.4) in
   $\bigoplus\Mat(\Lambda^j V)$,
   and  limit (3.6) in $\Mat (H_\nu)$ remain the same.

   \smallskip

 b) If we change the variable $z$ by the formula
   $$z\mapsto z^p,$$
   then the exponents $m_j$ are replaced by $pm_j$.
   All other data
   ( i.e., $\cP^\gamma$ and $\cL^j$)
   remain the same as above.

 c) Let us  change the variable $z$ by
   $$c\mapsto cz,\qquad\text{where}\quad c\in\C^*.$$
   Then $m_j$ do not change. The hinge $\cP^\gamma$
   transforms by the rule
   $$(P_1,\dots,P_\tau)\rightarrowtail
     (c^{-k_1} P_1,\dots, c^{-k_\tau} P_\tau).$$

   The collection
   $\cL(\cP^\gamma)=(\cL^0(\cP^\gamma),
   \dots,\cL^n(\cP^\gamma))\in \widetilde\Hinge_n$
    transforms by the rule
   $$(\cL^0,\cL^1,\cL^2,\dots)\rightarrowtail
   (c^0\cdot \cL^0, c^{-m_1}\cL^1,c^{-m_1-m_2}\cL^2,\dots).$$

   The operators
  $ [\rho_\nu(\gamma)]$ transform by the rule
     $$
     [\rho_\nu(\gamma)]
         \rightarrowtail c^{-\sum m_j \nu_j} [\rho_\nu(\gamma)]
     .$$

 {\bf 3.7. Sea urchin}.
  We define  the {\it sea urchin} $\SU_n$ as the union
  of $\GL_n(\C)$ and all the spikes $\sp[{\bold m}]$.

  The {\it spikes} $\sp[{\bold m}]$ are enumerated by
  the collections of integers
  $$ {\bold m}: m_1\ge m_2\ge\dots\ge m_n$$
  such that $m_1,\dots,m_n$ have no common divisor.
  Let $k_j$, $\alpha_j$ be the associated
  numbers (see 3.2). Points of the spike $\sp[{\bold m}]$
  are hinges  $\cP\in\cO[\alpha_1,\dots,\alpha_\tau]$
  defined up to the equivalence
     \begin{equation}(P_1,\dots,P_\tau)
 \sim    (c^{k_1} P_1,\dots, c^{k_\tau} P_\tau);\qquad
   c\in \C^*
   .\end{equation}

 {\sc Remark.} $\dim\GL_n=n^2$, and the dimension
 of all spikes is $n^2-1$.     \konets

 \smallskip

 {\bf 3.8. Limits of meromorphic curves in sea urchin.}
 Let $\gamma$ be a meromorphic curve in $\GL_n(\C)$.
 Let $\hat m_1,\dots, \hat m_n$ be its exponents,
 let $u$ be the greatest (positive) common divisor
 of the numbers $\hat m_j$. Then we define the collection
 $\bold m^\gamma$  by
 $$m_j=m_j^\gamma=\hat m_j/u$$
 The limit $\cP^\gamma$ of the curve $\gamma$ in the spike
 $\sp[{\bold m}]$ is defined by Proposition 3.2.

 By 3.6, the collection $\bold m^\gamma$ and  the limit
 of the curve $\gamma$ in $\sp[{\bold m}]$ do not depend on
  the parametrization
 of the curve $\gamma$.

 \smallskip

 {\bf 3.9. Projective compactifications of $\GL_n(\C)$.}
 Let $\rho_{\nu^{(1)}}$, \dots, $\rho_{\nu^{(\sigma)}}$
 be an irreducible holomorphic representations of $\GL_n(\C)$,
 where
 $$\nu^{(l)}:\,\nu_1^{(l)}\ge\dots\ge \nu_n^{(l)}$$
 are signatures (see 2.11).
 Let $H_{\nu^{(l)}}$ be the spaces
 of the representations
 $\rho_{\nu^{(l)}}$.
 Let
 $$\zeta:=\bigoplus_{l=1}^\sigma  \rho_{\nu^{(l)}};
   \qquad Z=\bigoplus_{l=1}^\sigma H_{\nu^{(l)}}.$$

 We have the maps
 \begin{equation}
  \GL_n(\C)\stackrel{\zeta}{\to}
      \Mat( Z)\setminus 0\to\PMat(Z).
 \end{equation}

Denote by $[\GL_n]_\zeta$ the closure of
the image of
$\GL_n(\C)$ in $\PMat(Z)$.
The spaces $\PMat(Z)$  are called
{\it projective compactifications}
of $\GL_n(\C)$.

\smallskip

{\bf 3.9. The canonical maps
$\pi_\zeta:\SU_n\to\PMat(Z)$.}
We must extend the map (3.8) to the spikes
$\sp[\bold m]$ of the sea urchin $\SU_n$.

Let us fix a spike   $\sp[\bold m]$. Consider
the numbers
$$v^{(l)}:= \sum_j m_j \nu_j^{(l)}.$$
Let
$$v:=\max_l v^{(l)}.$$
For a hinge $\cP\in  \sp[\bold m]$, we define the operator
$$
\zeta({\bold m};\cP)=\bigoplus B_l({\bold m},\cP):
    \bigoplus  H_{\nu^{(l)}}     \to   \bigoplus  H_{\nu^{(l)}}
,$$
where the operators
$$B_l=B_l({\bold m};\cP):
   H_{\nu^{(l)}}\to H_{\nu^{(l)}}$$ are given
by
$$B_l=\begin{cases}
         \rho_{\nu^{(l)}}(\cL(\cP^\gamma)), &\text{if}\quad v^{(l)}=v\\
                   0,   &\text{if}\quad v^{(l)} < v
        \end{cases}
        $$

 {\sc Remark.} The hinge $\cP\in \sp[{\bold m}]$ is defined up to the
equivalence
 (3.7). Hence the operators  $[\rho_{\nu^{(l)}}(\cP^\gamma)]$
 are defined up to the factors
 $c^{\sum_i m_j \nu_j^{(l)}}=c^{v^{(l)}}$.
 But for all nonzero $B_l$, these factors coincide,
 and hence $\zeta({\bold m};\cP)$ is a nonzero operator
  defined up
 to a factor $c^v$.        \konets

 \smallskip

 {\sc Theorem 3.5}  {\it The image of the sea urchin under
 the map $({\bold m},\cP)\mapsto \zeta({\bold m},\cP)$
 coincides with $[\GL_n]_\zeta$.}

 \smallskip

 {\sc Lemma 3.6.} {\it Let $z\mapsto \gamma(z)$ be a meromorphic
 family. Then the limit of $\zeta(\gamma(z))$ in
 $[\GL_n]_\zeta$ exists and coincides with
  $\zeta({\bold m}^\gamma, \cP^\gamma)$.}

 \smallskip

 {\sc Proof of Lemma 3.6.}
 $$\ord \rho_{\nu^{(l)}}\bigl(\gamma(z)\bigr)=
  \sum_i\nu_i^{(l)}m_i=v^{(l)}$$
  Thus our limit
  is
  $$\lim\limits_{z\to 0}z^{v}\zeta\bigl(\gamma(z)\bigr)
  =
  \lim\limits_{z\to 0}z^{v}\bigoplus_l \rho_{\nu^{(l)}}
  \bigl(\gamma(z)\bigr)$$
  and Lemma 3.4 implies the required result. \konets

 \smallskip

 {\sc Proof of Theorem 3.5.} The set
 $[\GL_n]_\zeta\setminus\GL_n(\C)$
 is a subvariety  of the projective variety
 $[\GL_n]_\zeta$. Thus any point
  of this set can be achieved
 along an algebraic curve $\gamma(z)\subset \GL_n(\C)$.
 By  Lemma 3.6,  this point is contained in the image of
 the sea urchin.

 Conversely, each point   $({\bold m}, \cP)$
  of the sea urchin can be achieved
 along an algebraic curve $\gamma(z)$.
 By Lemma 3.6, the image of the point
   $({\bold m}, \cP)$ can be achieved
   along  the curve $\zeta(\gamma(z))$
   and thus  $\zeta({\bold m}, \cP)\in[\GL_n]_\zeta $. \konets

 \smallskip

 {\sc Remark.} Let us give a description
  of the sea urchin
 in formal terms. The spikes are homogeneous
  spaces $\GL_n(\C)\times\GL_n(\C)/G_{\bold m}$,
  where the family of subgroups    $G_{\bold m}$
  can be obtained by the following algorithm.

  Consider two opposite parabolic subgroups
  $P_+$ and $P_-$ in
  $\GL_n(\C)$. Denote by $N_+$ and $N_-$
  the unipotent radicals in $P_+$, $P_-$.
  Denote by $Q$ the Levi subgroup ($Q\subset P_+$,
  $Q\subset P_-$).  Denote by $Z$
  the center of the Levi subgroup.
Let $S$ be an one-parametric subgroup in
$Z$. The say that $S$ is positive if
all root lying in $N_+$ are positive on
on the  generator of $S$.
Now we consider the subgroup
$$G^S\subset \GL_n(\C)\times\GL_n(\C)$$
consisting of all pairs
$$(g_1,g_2)=(n_- q_1, q_2 n_+)$$
where
$n_-\in N_-$, $n_+\in N_+$, $q_1, q_2\in Q$
and $q_1^{-1} q_2\in S$. Then
the family of the subgroups
$G^S$ coincides with the family $G_{\bold m}$.

This remark allows to extends the sea urchin construction
to an arbitrary complex semisimple group.
For all classical groups the De Concini--Procesi
compactification can be described on the hinge
language (see \cite{Ner3}, \cite{Ner5})
and the sea urchins
also can easily be
 described in these terms.  \konets

{ \sc
 Address (spring 2001):

  Erwin Schr\"odinger Institute for Mathematical
          Physics

 Boltzmanngasse, 9, Wien 1020, Austria

 \smallskip

 Permanent address:
  Institute of Theoretical and Experimental Physics,

  Bolshaya Cheremushkinskaya, 25

  Moscow 117259

  Russia}

  \smallskip

 {\tt e-mail neretin@main.mccme.rssi.ru, neretin@gate.itep.ru}

 \end{document}